\documentclass{tpms-l}
\usepackage{mathrsfs}
\newtheorem{theorem}{Theorem}[section]

\theoremstyle{definition}
\newtheorem{definition}[theorem]{Definition}
\newtheorem{example}[theorem]{Example}

\theoremstyle{remark}
\newtheorem{remark}[theorem]{Remark}

\numberwithin{equation}{section}

\begin{document}

\title[Recurrence of Markov chains]{A new criterion for recurrence of Markov chains with an infinitely countable set of states}

\author{Vyacheslav M. Abramov}
\address{24 Sagan Drive, Cranbourne North, Melbourne, Victoria-3977, Australia}
\email{vabramov126@gmail.com}



\date{10/APR/2024}



\begin{abstract}
For a class of irreducible  Markov chains with an infinitely countable set of states, we establish a new verifiable necessary and sufficient condition for recurrence and transience. We show that if one of the basic assumptions is not satisfied, then the statement of the theorem becomes invalid.
\end{abstract}

\maketitle


\section{Introduction}

\subsection{Formulation of the main result}
Let $\mathscr{P}$ be an irreducible Markov chain with the infinite set of states $\{0, 1, \ldots\}$. Let $p_{i,j}^{(n)}$ denote $n$-step transition probabilities from state $i$ to state $j$, $p_{i,j}=p_{i,j}^{(1)}$. Recall that a Markov chain is recurrent, if its initial state $0$ is visited infinitely many times with probability 1. It is well-known that the Markov chain $\mathscr{P}$ is recurrent, if and only if $\sum_{n=1}^{\infty}p_{0,0}^{(n)}=\infty$ (see \cite{Ch, N}). The application of this last condition is very limited, since the derivation of $p_{0,0}^{(n)}$ is generally combinatorially hard. The known applications are chiefly  related to simple random walks and can be found in many textbooks on the theory of Markov processes. The more extended study for the recurrence of Markov chains can be found in \cite[pp. 65, 66]{Ch}, where the aforementioned elementary result has been extended. Specifically, it was derived the system of equations that clarified the solution of this problem under quite general conditions. However, the derivation of explicit conditions for recurrence or transience from that system of equation is hard and can be provided for only the particular cases such as the Markov chains of the birth-and-death type, for which the aforementioned system of equations has been further applied in \cite[pp. 67--71]{Ch}.

In this paper, we suggest a new result that enables us to verify the condition of recurrence or transience for a relatively wide class of irreducible Markov chains. It turns out that the result presented here is an extension of the well-known one from the theory of birth-and-death processes. Furthermore, for a wide class of the Markov chains, the conditions for recurrence or transience are given in the same terms as those for the birth-and-death chains. In other words, we prove that a Markov chain of the considered class is recurrent if and only if a specifically defined birth-and-death chain is recurrent.

\begin{definition}
 The states of Markov chain $\mathscr{P}$ are said to form connected
 domain if for any $i \geq 1$ there are $j_1(i)$ and $j_2(i)$, $0 \leq j_1(i) < i < j_2(i) \leq\infty$
 such that $p_{i,k} > 0$ for all $j_1(i) \leq k \leq i - 1$ and $i + 1 \leq k \leq j_2(i)$, and $\sum_{k=j_1(i)}^{j_2(i)}p_{i,k}=1$.
\end{definition}

A trivial example of a Markov chain, the states of which form connected domain is the birth-and-death type chain, where the diagonal elements of the transition matrix are not necessarily zeros.
Another example of a Markov chain, the states of which form connected domain
 is a Markov chain, where $p_{i,j} > 0$ for all $i,j \geq 1$ with $i\neq j$. The examples
 of Markov chains the states of which do not form connected domain are
 discussed in Section \ref{S3}.

 \begin{theorem}\label{Th1}
  Assume that the states of the Markov chain $\mathscr{P}$ form connected domain, and
  $\sum_{j=1}^{\infty}jp_{0,j}<\infty$.
  Denote $e_{i}^-=\sum_{j=1}^{\infty}jp_{i+j-1,i-1}$ and $e_{i}^+=\sum_{j=1}^{i+1}jp_{i-j+1,i+1}$,
  $i\geq1$, assuming that $e_{i}^-<\infty$.
  Assume also that $1-p_{i,i}>\epsilon$ for all $i\geq0$, where $\epsilon$ is some positive value.
  Then, the Markov chain $\mathscr{P}$ is recurrent if and only if
  \[
  \sum_{n=1}^\infty\prod_{i=1}^{n}\frac{e_i^-}{e_i^+}=\infty.
  \]
\end{theorem}

\begin{remark}
The assumptions $\sum_{j=1}^{\infty}jp_{0,j}<\infty$ and $e_{i}^-<\infty$, $i\geq1$,
are technical, since the infinite series $\sum_{j=1}^{\infty}jp_{0,j}$ and $\sum_{j=1}^{\infty}jp_{i+j-1,i-1}$
appear explicitly in the proof of the theorem. It is possible to avoid these assumptions using additional arguments. Specifically, one can consider the series of Markov chains, in which the required parameters are assumed to increase to infinity and arrive at the conclusion in the limit case.
\end{remark}

The well-known corollary from Theorem \ref{Th1} is the condition of recurrence or transience for birth-and-death chains (e.g., \cite[pp. 66--71]{Ch}).

\subsection{Methodology of the study} The direct methods for the analysis of Markov chains, such as suggested in \cite[pp, 65, 66]{Ch} for instance, seem to be insufficient for the solution of the problem considered in the present paper. The system of equations obtained in \cite[p. 66]{Ch} can be served as a starting point of the study that may require tedious analytical derivations and finally lead to the equations similar to these being used in the paper.

We decided to follow another way. The analysis of this paper is based on the methods of stochastic calculus and the theory of martingales in continuous time and involves the theory of point processes. For this purpose, the original Markov chains is extended to the Markov processes in continuous time. Such approach has an advantage, since the results for the birth-and-death processes and point processes are widely presented in the literature and widely known, while the similar study based on the discrete time processes seems do not yield any simplifications.

For those continuous time Markov processes, we derive the system of stochastic differential equations. We obtain the equations that describe limit behavior of the process and prove that the conditions for recurrence or transience are the same as those for the Markov processes of the birth-and-death type.

\subsection{Outline of the paper}
The rest of this paper is organized as follows. In Section \ref{S2}, we prove the theorem. Our methodology is first demonstrated on elementary Example \ref{ex1} and then developed for a general type Example \ref{ex3} that closely describes the proof in the general case. In Section \ref{S3}, we provide simple but nontrivial examples of recurrent and transient Markov chains that demonstrate how the statement of the theorem works.
In Section \ref{S4}, we provide counterexamples that demonstrate that the key assumption of the theorem stating that a Markov chain forms a connected domain, is important. In Section \ref{S5}, we discuss our findings.

\section{Proof of the theorem}\label{S2} We start the proof from the preliminary information given in Sections \ref{S2.1} and \ref{S2.2}. In section \ref{S2.1}, we recall information about  Markov chains, the transition probability matrix of which describes a birth-and-death process. In Section \ref{S2.2}, some natural extension of those Markov chains helping us to further simplify the proof is provided. In Section \ref{S2.3} the main part of the proof is given. We provide a number of illustrative examples that clarify the further steps of the proof.

\subsection{Birth-and-death process}\label{S2.1} The  important class of Markov chains we deal with in this paper is the Markov chain $\mathscr{X}$ that is described by the following transition probability matrix:
\begin{equation}\label{2}
\left(\begin{matrix}0 &1 &0 &0 &0 &0 &\ldots\\
q_1 &0 &p_1 &0 &0 &0 &\ldots\\
0 &q_2 &0 &p_2 &0 &0 &\ldots\\
0 &0 &q_3 &0 &p_3  &0 &\ldots\\
\vdots &\vdots &\vdots &\ddots &\ddots &\ddots &\vdots \end{matrix}\right), \quad q_i+p_i=1, \quad i\geq1.
\end{equation}
It describes transition probabilities for a birth-and-death process. It is well-known that the Markov chain $\mathscr{X}$ is recurrent if and only if
$\sum_{n=1}^{\infty}\prod_{i=1}^{n}q_i/p_i=\infty$ (see \cite{KM}, \cite[p. 67--71]{Ch}, \cite[p. 286]{D}).

\subsection{Simplification of the problem}\label{S2.2}
We start from the following notice related to the natural extension of \eqref{2}. Let $c_1$, $c_2$,\ldots be a sequence, $0<\epsilon<c_i\leq1$, $i\geq1$. Consider the Markov chain $\mathscr{Y}$ with the transition probabilities matrix
\begin{equation}\label{4}
\left(\begin{matrix}0 &1 &0 &0 &0 &0 &\ldots\\
c_1q_1 &1-c_1 &c_1p_1 &0 &0 &0 &\ldots\\
0 &c_2q_2 &1-c_2 &c_2p_2 &0 &0 &\ldots\\
0 &0 &c_3q_3 &1-c_3 &c_3p_3  &0 &\ldots\\
\vdots &\vdots &\vdots &\ddots &\ddots &\ddots &\vdots \end{matrix}\right).
\end{equation}

We prove first that the Markov chain $\mathscr{Y}$ is recurrent if and only if the Markov chain $\mathscr{X}$ is recurrent. Let $x_{i,j}$ denote the transition probability from state $i$ to state $j$ of the Markov chain $\mathscr{X}$, and let $y_{i,j}$ denote the transition probability from state $i$ to state $j$ of the Markov chain $\mathscr{Y}$. Apparently, that $x_{i,i-1}$ and $x_{i,i+1}$ can be considered as the conditional transition probabilities in the Markov chain $\mathscr{Y}$ given that the state $i$ of this Markov chain changes. Therefore, if the initial state of the Markov chain $\mathscr{X}$ is visited only finite number of times with probability 1, the same is true for the Markov chain $\mathscr{Y}$. Furthermore, the number of visits state $0$ for both of the Markov chains $\mathscr{X}$ and $\mathscr{Y}$ has the same distribution.
Otherwise, if the initial state of the Markov chain $\mathscr{X}$ is visited infinitely many number of times with probability 1, then the same is true for the Markov chain $\mathscr{Y}$ and vice versa.

Based on this specific example, one can assume without loss of generality that $p_{i,i}=0$ for all $i$, and under this additional assumption we shall prove our main result.

\subsection{Stochastic Chapman-Kolmogorov equations}\label{S2.3} The variant of the Chapman-Kolmogorov equations used in this section is a stochastic version of these equations that was originally used in \cite{KL} for a specified queueing network of a closed type.

We show that the structure of the equations for the steady state probabilities can be reduced to the same structure as that of the birth-and-death processes.
For this purpose, it is convenient to switch from the original Markov chain to the continuous time Markov process as follows. We assume that the times between transitions all are independent and exponentially distributed with mean $1/2$.

First, we give three examples to make our further arguments clearer. The first two examples describe elementary models. The third one describes a natural model and is closer to the formulated result. The further elementary examples of Section \ref{S3} are based on our general result.

\begin{example}\label{ex1}
Let $\mathscr{A}$ be a birth-and-death process, where the birth and death rates are $\lambda_n>0$ and $\mu_n>0$ when there are $n$ particles in the population. 
Let $\mathscr{B}$ be another stochastic process with the following behavior. When there are $n$ particles in the population, there can occur either two births simultaneously (i.e. the birth of twins) or one death. The birth rate is assumed to be $\lambda_{n}/2$ (i.e. times between births of twins are exponentially distributed with mean $2/\lambda_n$) while the death rate is $\mu_n$. Let $\lambda_0$ be an arbitrary positive regeneration parameter for the birth-and-death process $\mathscr{A}$ and let $\lambda_0/2$ be such one for the process $\mathscr{B}$, when there are no individuals in these populations. According to the traditional definition, the stochastic process $\mathscr{B}$ is not a birth-and-death process. However, for our convenience it will be further called birth-and-death process as well. It will be shown in this example that the birth-and-death process $\mathscr{B}$ is recurrent if and only if the birth-and-death process $\mathscr{A}$ is recurrent.

Note that the models $\mathscr{A}$ and $\mathscr{B}$ are described in the terms of Markov chain transition probability as follows. The model $\mathscr{A}$ is described by matrix \eqref{2}, where $p_i=\lambda_i/(\lambda_i+\mu_i)$, $q_i=\mu_i/(\lambda_i+\mu_i)$. Then for the model $\mathscr{B}$ the matrix of transition probability is
\begin{equation}\label{10}
\left(\begin{matrix}0 &1 &0 &0 &0 &0 &0&\ldots\\
c_1q_1 &0 &0 &\frac{c_1p_1}{2} &0 &0 &0&\ldots\\
0 &c_2q_2 &0 &0 &\frac{c_2p_2}{2} &0 &0 &\ldots\\
0 &0 &c_3q_3 &0 &0  &\frac{c_3p_3}{2} &0 &\ldots\\
\vdots &\vdots &\vdots &\ddots &\ddots &\ddots &\ddots &\vdots \end{matrix}\right),
\end{equation}
with the same values $p_i$ and $q_i$ and with the normalization constants $c_i$ such that the sums in each row are equal to $1$. Note that $c_iq_i=2\mu_i/(\lambda_i+2\mu_i)$, and $c_ip_i/2=\lambda_i/(\lambda_i+2\mu_i)$, $i\geq1$.

 The system of the equations for the birth-and-death processes that describes the process $\mathscr{A}$ is
\begin{eqnarray}
P_{n}\left(\lambda_n+\mu_n\right)&=&\mu_{n+1}P_{n+1}+\lambda_{n-1}P_{n-1}, \ n\geq1,\label{e1}\\
\lambda_0P_0&=&\mu_1P_1,\label{e1.1}
\end{eqnarray}
where $P_0, P_1,\ldots$ denote the steady state probabilities, if they exist (e.g., \cite[p. 285]{D}) . Existence of the steady state probabilities means that $\sum_{n=0}^\infty P_n=1$. Otherwise, $P_n=0$ for all $n\geq0$.

For the process $\mathscr{B}$, we will derive the system of equations
\begin{eqnarray}
Q_{n}(\lambda_n+\mu_n)&=&\mu_{n+1}Q_{n+1}+\lambda_{n-2}Q_{n-2},\ n\geq2,\label{e2}\\
Q_1(\lambda_1+\mu_1)&=&\mu_2Q_2,\label{e2.1}\\
\lambda_0Q_0&=&\mu_1Q_1,\label{e2.2}
\end{eqnarray}
where $Q_0, Q_1,\ldots$ denote the steady state probabilities, if they exist.  Note that according to \eqref{e2}--\eqref{e2.2} if the steady state probabilities exist, then $\sum_{n=0}^\infty Q_n=1$, otherwise $Q_n=0$ for all $n\geq0$.

There are several ways to solve the system of equations \eqref{e2}--\eqref{e2.2}. One of the known methods is based on the system of difference equations applied to the Markov chain transition probabilities given by \eqref{10}. Another equivalent method related to the continuous time process is based on the system of Chapman-Kolmogorov equations. However the justification of \eqref{e2}--\eqref{e2.2} needs special steps, since direct application of the same definition as in the case of the birth-and-death equations would require to justify the presence of a factor $2$ before every $\mu_n$.
In our considerations, the system of equations \eqref{e2}--\eqref{e2.2} is considered as particular of the more general result that will be provided later based on \textit{stochastic system of Chapman-Kolmogorov equations}. Application of the stochastic system of Chapman-Kolmogorov equations is more profitable compared to the traditional one for the following reason. 
The stochastic system of Chapman-Kolmogorov equations is written for $\mathsf{I}\{Z(t)=n\}$ of the continuous Markov process $Z(t)$ that describes the number of individuals in the population at time $t$. Taking expectation from the both parts of each equation will lead to the usual system of Chapman-Kolmogorov equations. Then, the advantage of the application of stochastic system of Chapman-Kolmogorov equations compared to the traditional one is that one can take expectation not immediately, but in our convenience after deriving necessary representations. Specifically, we first use the semimartingale decomposition for $\mathrm{I}\{Z(t)=n\}$, and then, when we obtain a more convenient representation, we can apply the expectation for the both sides of the equations. In this case we do not need to solve the system of differential equations or prove its convergence to the linear system of equations as it is required under the traditional approach, but arrive at the necessary linear system of equations directly.

Let us now show that the birth-and-death process $\mathscr{B}$ is recurrent if and only if the same is true for the birth-and-death process $\mathscr{A}$.

Indeed, together with \eqref{e1}, \eqref{e1.1} we also have
\begin{equation}\label{1}
\lambda_0P_0+\sum_{n=1}^{\infty}P_n(\lambda_n+\mu_n)=\mu_1P_1+\sum_{n=1}^{\infty}(\lambda_{n-1}P_{n-1}+\mu_{n+1}P_{n+1}),
\end{equation}
and together with \eqref{e2}--\eqref{e2.2} we also have

\begin{align}
\lambda_0Q_0+\sum_{n=1}^{\infty}Q_n(\lambda_n+\mu_n)&=\mu_1Q_1+\mu_2Q_2+\sum_{n=2}^{\infty}(\lambda_{n-2}Q_{n-2}+\mu_{n+1}Q_{n+1})\nonumber\\
&=\mu_1Q_1+\sum_{n=1}^{\infty}(\lambda_{n-1}Q_{n-1}+\mu_{n+1}Q_{n+1}).\label{3}
\end{align}

So, we arrived at identical equations \eqref{1} and \eqref{3}. Hence, $\sum_{n=0}^{\infty}Q_n=1$ if and only if $\sum_{n=0}^{\infty}P_n=1$. (The last is true, since $\sum_{n=0}^{\infty}Q_n$ can take the only values $0$ or $1$.)
\end{example}

Apparently, the model for the birth-and-death process $\mathscr{B}$ considered in Example \ref{ex1}
can be further extended (for instance we can assume that the birth rate is $\lambda_n/l$ and any new birth is given for $l$ new particles, while the death rate $\mu_n$ remains the same as above), and from the point of view of recurrence or transience we will have the class of equivalent processes satisfying the same property.

Consider another example close to Example \ref{ex1}.

\begin{example}\label{ex2}
Let $\mathscr{C}$ be a stochastic process of the birth-and-death type with parameters of birth $\lambda_n/2$ and death $\mu_n$. While the deaths in state $n$ occur in standard way and the death rate is $\mu_n$, the births  occur according to the following scheme. It is either one birth with probability half or three birth with the same probability half. It follows from the general result provided in this paper that the stochastic process $\mathscr{C}$ is recurrent if and only if the birth-and-death process $\mathscr{A}$ is recurrent. That is, the $\mathscr{A}$, $\mathscr{B}$ and $\mathscr{C}$ all are either recurrent or transient.

The steady state probabilities of the process $\mathscr{C}$ are described by the following system of equations:
Here we have the following system of difference equations
\begin{eqnarray}
Q_{n}\left(\lambda_n+\mu_n\right)&=&\mu_{n+1}Q_{n+1}+\frac{\lambda_{n-1}}{4}Q_{n-1}+\frac{3\lambda_{n-3}}{4}Q_{n-3},\ n\geq3,\label{e3}\\
Q_2(\lambda_2+\mu_2)&=&\mu_3Q_3+\frac{\lambda_1}{4}Q_1,\label{e3.1}\\
Q_1(\lambda_1+\mu_1)&=&\mu_2Q_2+\frac{\lambda_0}{4}Q_0,\label{e3.2}\\
\lambda_0Q_0&=&\mu_1Q_1\label{e3.3}
\end{eqnarray}
where $Q_0$, $Q_1$,\ldots denote the steady state probabilities. The system of equations \eqref{e3}--\eqref{e3.3} follows from the derivation presented later in the proof for a more general model.
\end{example}

With $p_i=\lambda_i/(\lambda_i+\mu_i)$, $q_i=\lambda_i/(\lambda_i+\mu_i)$ the matrix of transition probabilities for the model $\mathscr{C}$ is
\begin{equation*}
\left(\begin{matrix}0 &1 &0 &0 &0 &0 &0&\ldots\\
c_1q_1 &0 &\frac{c_1p_1}{4} &0 &\frac{c_1p_1}{4} &0 &0&\ldots\\
0 &c_2q_2 &0 &\frac{c_2p_2}{4} &0 &\frac{c_2p_2}{4} &0 &\ldots\\
\vdots &\vdots &\ddots &\ddots &\ddots &\ddots &\ddots &\vdots \end{matrix}\right).
\end{equation*}
where the constants $c_1$, $c_2$,\ldots are the same as for matrix \eqref{10}. That is, $c_iq_i=2\mu_i/(\lambda_i+2\mu_i)$ and $c_ip_i/4=\lambda_i/(2\lambda_i+4\mu_i)$, $i\geq1$.

\smallskip

In the following study, the class of equivalent processes should strongly contain the models of Markov chains that form connected domain only, although the provided proof for the birth-and-death process $\mathscr{B}$ did not use any restriction.  Despite the processes $\mathscr{B}$ and $\mathscr{C}$ do not form connected domain, the proof provided later for models that form connected domain will remain true for these processes. This is because their structure is simple.

\begin{example}\label{ex3}
We consider a typical example of a Markov chain that forms connected domain (call it $\mathscr{D}$), the matrix of transition probabilities of which is
\begin{equation*}
\left(\begin{matrix}0 &1 &0 &0 &0 &0   &\ldots\\
q_{1,1} &0 &p_{1,1} &p_{1,2} &p_{1,3} &p_{1,4}   &\ldots\\
q_{2,1} &q_{2,2} &0 &p_{2,1} &p_{2,2} &p_{2,3}   &\ldots\\
q_{3,1} &q_{3,2} &q_{3,3} &0 &p_{3,1}  &p_{3,2}  &\ldots\\
\vdots &\vdots &\vdots &\vdots &\vdots  &\vdots  &\ddots \end{matrix}\right),
\end{equation*}
where all $p_{i,j}$ and $q_{i,j}$ that appear in the matrix are assumed to be strictly positive and $\sum_{j=1}^{i}q_{i,j}+\sum_{j=1}^{\infty}p_{i,j}=1$, $i=1,2,\ldots$,
and the Markov chain of the birth-and-death type (call it $\mathscr{E}$), the matrix of transition probabilities of which is
\begin{equation*}
\left(\begin{matrix}0 &1 &0 &0 &0&\ldots\\
r_{1} &0 &s_1 &0  &0&\ldots\\
0 &r_2 &0 &s_2 &0&\ldots\\
\vdots &\vdots &\ddots &\ddots &\ddots &\vdots\end{matrix}\right),
\end{equation*}
where
\[
r_n=c_n\sum_{i=n}^{\infty}(i-n+1)q_{i,n}, \quad s_n=c_n\sum_{i=1}^{n}(n-i+1)p_{i,n-i+1}, \quad n\geq1,
\]
where the normalization constants $c_i$ are given such that sums $r_i+s_i$ in each row $i\geq1$ are equal to $1$.

Our aim is to demonstrate that both Markov chains $\mathscr{D}$ and $\mathscr{E}$ are either recurrent or transient.

From the theory of birth-and-death processes, for the Markov chain $\mathscr{E}$ we have the system of equations:
\begin{equation}\label{e12}
\begin{aligned}
&\left(\sum_{i=n}^{\infty}(i-n+1)q_{i,n}+\sum_{i=1}^{n}(n-i+1)p_{i,n-i+1}\right)P_n\\
&=P_{n-1}\sum_{i=0}^{n-1}(n-i)p_{i,n-i}+P_{n+1}\sum_{i=n+1}^{\infty}(i-n)q_{i,n+1},\quad n\geq1,
\end{aligned}
\end{equation}
\begin{equation}\label{e12.1}
P_0=P_1\sum_{i=1}^{\infty}iq_{i,1},
\end{equation}
where $P_0, P_1,\ldots$ denote the steady state probabilities (if exist).

We shall prove below that, for the Markov chain $\mathscr{D}$ the following system of equation is satisfied
\begin{equation}\label{e13}
\begin{aligned}
&\left(\sum_{i=n}^{\infty}(i-n+1)q_{i,n}+\sum_{i=1}^{n}(n-i+1)p_{i,n-i+1}\right)Q_n\\
&=\sum_{i=0}^{n-1}(n-i)p_{i,n-i}Q_{i}+\sum_{i=n+1}^{\infty}(i-n)q_{i,n+1}Q_{i},\quad n\geq1,
\end{aligned}
\end{equation}
\begin{equation}\label{e13.1}
Q_0=Q_1\sum_{i=1}^{\infty}iq_{i,1},
\end{equation}
where $Q_0, Q_1,\ldots$ denote the steady state probabilities (if exist).

The system of equations that has been derived in \cite[p. 66]{Ch} for limit probabilities seems do not make possible to arrive in some way at \eqref{e13}, \eqref{e13.1}. As well, the system of equations cannot be derived directly by the same method as it is given in Example \ref{ex1}. The derivation that will be provided below involves the theory of point processes and requires special measurability arguments for the martingales that will appear in the derivation. Therefore our approach will be based on the extension of the method that has been suggested in the derivation of \eqref{e2}--\eqref{e2.2}.

\smallskip
\textit{Proof of \eqref{e13}, \eqref{e13.1}}. As before, we extend the Markov chain to continuous time Markov process, assuming that the times between the jumps are exponentially distributed with parameter $1/2$. Denote this process by $Z(\tau)$ assuming that $Z(0)=0$. For any time $\tau>0$ we set $Z(\tau)=\Lambda(\tau)-\mathrm{M}(\tau)$, where $\Lambda(\tau)$ and $\mathrm{M}(\tau)$ are state-dependent compound Poisson processes depending on the state of the process at time $\tau$. In the sequel, we will introduce a series of counting processes that will be denoted by $\Lambda_{n}(\tau)$ and $\mathrm{M}_{n}(\tau)$, $n\geq0$, associated with $\Lambda(\tau)$ and $\mathrm{M}(\tau)$, respectively.

Let $\tau_1$, $\tau_2$,\ldots denote the times of jumps when $\mathsf{I}\{Z(\tau_i-)=n, Z(\tau_i)>n\}=1$. These times of jumps will be denoted by $\tau_1^{(n)}$, $\tau_2^{(n)}$,\ldots.
Similarly, let $\eta_1$, $\eta_2$,\ldots denote the times of jumps when $\mathsf{I}\{Z(\eta_i-)=n, Z(\eta_i)<n\}=1$. These times of jumps will be denoted by $\eta_1^{(n)}$, $\eta_2^{(n)}$,\ldots.
Denote by $\triangle Z(\tau_i^{(n)})=Z(\tau_i^{(n)})-Z(\tau_i^{(n)}-)$ the jumps at times $\tau_i^{(n)}$, $i=1,2,\ldots$ and by $\triangle Z(\eta_i^{(n)})=Z(\eta_i^{(n)})-Z(\eta_i^{(n)}-)$ the jumps at times $\eta_i^{(n)}$, $i=1,2,\ldots$. The sequence of $\triangle Z(\tau_i^{(n)})$ (containing positive values) is denoted by $\triangle \Lambda_n(\tau_i^{(n)})$ and associated with the process $\Lambda_n(\tau)$, where the process $\Lambda_n(\tau)$ itself is obtained by deleting all the intervals of the process $\Lambda(\tau)$ when $\mathsf{I}\{Z(\tau-)\neq n\}$ and merging the ends, with following re-scaling times.
The sequence $\triangle Z(\eta_i^{(n)})$ (containing negative values) taken with the opposite (positive) sign is denoted by $\triangle \mathrm{M}_n(\eta_i^{(n)})$ and associated with the process $\mathrm{M}_n(\tau)$, where the process $\mathrm{M}_n(\tau)$ is obtained similarly to that of $\Lambda_n(\tau)$. It makes sense to renumber the elements of sequences $\triangle \Lambda_n(\tau_i^{(n)})$ and $\triangle \mathrm{M}_n(\eta_i^{(n)})$ and denote them by $X_1^{(n)}$, $X_2^{(n)}$,\ldots, and $Y_1^{(n)}$, $Y_2^{(n)}$,\ldots, respectively. Each of the sequences $\{X_i^{(n)}\}_{i\geq1}$, $\{Y_i^{(n)}\}_{i\geq1}$ is assumed to consist of i.i.d. random variables.

We have:
\begin{eqnarray*}
\Lambda_n(\tau)&=&\sum_{i=1}^{N_n^\Lambda(\tau)}X_i^{(n)}=\int_{0}^{\tau}X_{N_n^\Lambda(t)}^{(n)}\mathrm{d}N_n^\Lambda(t),\\
\mathrm{M}_n(\tau)&=&\sum_{i=1}^{N_n^{\mathrm{M}}(\tau)}Y_i^{(n)}=\int_{0}^{\tau}Y_{N_n^{\mathrm{M}}(t)}^{(n)}\mathrm{d}N_n^{\mathrm{M}}(t),
\end{eqnarray*}
where $N_n^\Lambda(\tau)$, $N_n^{\mathrm{M}}(\tau)$ are the point processes associated with $\Lambda_n(\tau)$ and $\mathrm{M}_n(\tau)$, respectively.

The system of equations for the indicators $\mathsf{I}\{Z(\tau)=k\}$ with probability 1 is
\begin{equation}\label{12}
\begin{aligned}
\mathsf{I}\{Z(\tau)=n\}=&\sum_{j=n+1}^{\infty}\mathsf{I}\{Z(\tau-)=j\}\mathsf{I}\{\triangle\mathrm{M}_j(\tau)=j-n\}\\
&+\sum_{j=0}^{n-1}\mathsf{I}\{Z(\tau-)=j\}\mathsf{I}\{\triangle\Lambda_j(\tau)=n-j\}\\
&+\mathsf{I}\{Z(\tau-)=n\}\mathsf{I}\{\triangle\mathrm{M}_n(\tau)=0\}\mathsf{I}\{\triangle\Lambda_n(\tau)=0\}\\
=&I_1+I_2+I_3.
\end{aligned}
\end{equation}

\smallskip
Note that
\begin{align*}
\triangle\Lambda_n(\tau)&=X_{N_n^\Lambda(\tau)}^{(n)}\triangle N_n^\Lambda(\tau),\\
\triangle\mathrm{M}_n(\tau)&=Y_{N_n^\mathrm{M}(\tau)}^{(n)}\triangle N_n^\mathrm{M}(\tau),
\end{align*}
where $\triangle N_n^\Lambda(\tau)$ and $\triangle N_n^\mathrm{M}(\tau)$ denote the increments of the corresponding point processes taking value $1$ (a jump occurs exactly at time $\tau$ when the system in state $n$ at time $\tau-$) or $0$ (otherwise).

We have the following relations:
\begin{align}
\mathsf{I}\{\triangle\Lambda_j(\tau)=n-j\}&=\mathsf{I}\{X_{N_j^\Lambda(\tau)}^{(j)}=n-j\}\triangle N_j^\Lambda(\tau)\label{14}\\
&=\mathsf{I}\{X_{N_j^\Lambda(\tau-)+1}^{(j)}=n-j\}\triangle N_j^\Lambda(\tau),\label{16}
\end{align}
and
\begin{align}
\mathsf{I}\{\triangle\mathrm{M}_j(\tau)=j-n\}&=\mathsf{I}\{Y_{N_j^\mathrm{M}(\tau)}^{(j)}=j-n\}\triangle N_j^\mathrm{M}(\tau)\label{22}\\
&=\mathsf{I}\{Y_{N_j^\mathrm{M}(\tau-)+1}^{(j)}=j-n\}\triangle N_j^\mathrm{M}(\tau).\label{24}
\end{align}
Representations \eqref{14} and \eqref{16} are equivalent, since
\[
\mathsf{I}\{X_{N_j^\Lambda(\tau)}^{(j)}=n-j\}\triangle N_j^\Lambda(\tau)=\mathsf{I}\{X_{N_j^\Lambda(\tau-)+1}^{(j)}=n-j\}\triangle N_j^\Lambda(\tau)
\]
if $\triangle N_j^\Lambda(\tau)=1$. The same observation pertains to \eqref{22} and \eqref{24}.  The significance of
these relations will be seen later when we discuss the predictability of $\mathsf{I}\{X_{N_j^\Lambda(\tau-)+1}^{(j)}=n-j\}$ and $\mathsf{I}\{Y_{N_j^\mathrm{M}(\tau-)+1}^{(j)}=j-n\}$ for the martingale arguments used in the paper.
(For the further measurability arguments, we need to specify all the processes that
appear in the stochastic equations at time $\tau-$, that is, immediately before a possible
jump at time $\tau$.)

In addition to these equation, with probability $1$ we have
\begin{equation}\label{26}
\begin{aligned}
\mathsf{I}\{\Lambda_n(\tau)=0\}\mathsf{I}\{\mathrm{M}_n(\tau)=0\}=&[1-X_{N_n^\Lambda(\tau)}^{(n)}\triangle N_n^\Lambda(\tau)]\\
&\times[1-Y_{N_n^\mathrm{M}(\tau)}^{(n)}\triangle N_n^\mathrm{M}(\tau)].
\end{aligned}
\end{equation}

To derive stochastic differential equations, relations \eqref{16}, \eqref{24} and \eqref{26} are rewritten below in the form of stochastic differentials. From \eqref{16} and \eqref{24} we obtain:
\begin{equation}\label{28}
\mathrm{d}\mathsf{I}\{\triangle\Lambda_j(\tau)=n-j\}=\mathsf{I}\{X_{N_j^\Lambda(\tau-)+1}^{(j)}=n-j\}\mathrm{d} N_j^\Lambda(\tau),
\end{equation}
\begin{equation}\label{30}
\mathrm{d}\mathsf{I}\{\triangle\mathrm{M}_j(\tau)=j-n\}=\mathsf{I}\{Y_{N_j^\mathrm{M}(\tau-)+1}^{(j)}=j-n\}\mathrm{d} N_j^\mathrm{M}(\tau).
\end{equation}

In turn, equation \eqref{26} implies that
\begin{equation}\label{32}
\begin{aligned}
\mathrm{d}(\mathsf{I}\{\Lambda_n(\tau)=0\}\mathsf{I}\{\mathrm{M}_n(\tau)=0\})=&-X_{N_n^\Lambda(\tau)}^{(n)}\mathrm{d}N_n^\Lambda(\tau)
-Y_{N_n^\mathrm{M}(\tau)}^{(n)}\mathrm{d}N_n^\mathrm{M}(\tau)\\
=&-X_{N_n^\Lambda(\tau-)}^{(n)}\mathrm{d}N_n^\Lambda(\tau)
-Y_{N_n^\mathrm{M}(\tau-)}^{(n)}\mathrm{d}N_n^\mathrm{M}(\tau),
\end{aligned}
\end{equation}
since $X_{N_n^\Lambda(\tau-)}^{(n)}=X_{N_n^\Lambda(\tau)}^{(n)}$ and $Y_{N_n^\mathrm{M}(\tau-)}^{(n)}=Y_{N_n^\mathrm{M}(\tau)}^{(n)}$ in continuity point $\tau$.

According to representations \eqref{28}--\eqref{32}, with probability $1$ the differential form of equation \eqref{12} is
\begin{equation}\label{34}
\begin{aligned}
&\mathrm{d}\mathsf{I}\{Z(\tau)=n\}\\
=&\sum_{j=n+1}^{\infty}\mathsf{I}\{Z(\tau-)=j\}\mathsf{I}\{Y_{N_j^\mathrm{M}(\tau-)+1}^{(j)}=(j-n)\}\mathrm{d}N_j^M(\tau)\\
&+\sum_{j=0}^{n-1}\mathsf{I}\{Z(\tau-)=j\}\mathsf{I}\{X_{N_j^\Lambda(\tau-)+1}^{(j)}=(n-j)\}\mathrm{d}N_j^\Lambda(\tau)\\
&-\mathsf{I}\{Z(\tau-)=n\}[X_{N_n^\Lambda(\tau-)}^{(n)}\mathrm{d}N_n^\Lambda(\tau)+Y_{N_n^\mathrm{M}(\tau-)}^{(n)}\mathrm{d}N_n^\mathrm{M}(\tau)].
\end{aligned}
\end{equation}

In the sequel, we will need in the Doob-Meyer semimartingale decomposition (see \cite{JS, LS, P}) for the processes $N_j^\Lambda(\tau)$ and $N_j^\mathrm{M}(\tau)$ with respect to the filtration $\mathscr{F}_\tau$ defined as
\[
\mathscr{F}_\tau=\sigma\{Z(t), 0\leq t\leq\tau; X_i^{(n)}, Y_j^{(n)}, 1\leq i,j< \infty, n\geq 0\}.
\]
For this filtration, the sample paths of $Z(t)$, $N_i^\Lambda(t)$, $N_j^\mathrm{M}(t)$, $0\leq i,j<\infty$ are known up to time $\tau$, and all future jump sizes after time $\tau$ of $N_i^\Lambda(t)$, $N_j^\mathrm{M}(t)$, $0\leq i,j<\infty$ are known, but the locations of jumps are unknown.

The Doob-Meyer semimartingale decompositions for the Poisson processes $N_j^\Lambda(t)$, $N_j^\mathrm{M}(t)$, the rates of which is $1$, are
\[
N_j^\Lambda(\tau)=\tau+\boldsymbol{M}_{N_j^\Lambda(\tau)}, \quad N_j^\mathrm{M}(\tau)=\tau+\boldsymbol{M}_{N_j^\mathrm{M}(\tau)},
\]
where $\boldsymbol{M}_{N_j^\Lambda(\tau)}$ and $\boldsymbol{M}_{N_j^\mathrm{M}(\tau)}$ are square integrable martingales adapted to the filtration $\mathscr{F}_\tau$. Substitution of this for \eqref{34} yields:
\begin{equation}\label{36}
\begin{aligned}
&\mathrm{d}\mathsf{I}\{Z(\tau)=n\}\\
=&\sum_{j=n+1}^{\infty}\mathsf{I}\{Z(\tau-)=j\}\mathsf{I}\{Y_{N_j^\mathrm{M}(\tau-)+1}^{(j)}=(j-n)\}\mathrm{d}\tau\\
&+\sum_{j=0}^{n-1}\mathsf{I}\{Z(\tau-)=j\}\mathsf{I}\{X_{N_j^\Lambda(\tau-)+1}^{(j)}=(n-j)\}\mathrm{d}\tau\\
&-\mathsf{I}\{Z(\tau-)=n\}\left(X_{N_n^\Lambda(\tau-)}^{(n)}+Y_{N_n^\mathrm{M}(\tau-)}^{(n)}\right)\mathrm{d}\tau\\
&+\sum_{j=n+1}^{\infty}\mathsf{I}\{Z(\tau-)=j\}\mathsf{I}\{Y_{N_j^\mathrm{M}(\tau-)+1}^{(j)}=(j-n)\}\mathrm{d}\boldsymbol{M}_{N_j^M(\tau)}\\
&+\sum_{j=0}^{n-1}\mathsf{I}\{Z(\tau-)=j\}\mathsf{I}\{X_{N_j^\Lambda(\tau-)+1}^{(j)}=(n-j)\}\mathrm{d}\boldsymbol{M}_{N_j^\Lambda(\tau)}\\
&-\mathsf{I}\{Z(\tau-)=n\}[X_{N_n^\Lambda(\tau-)}^{(n)}\mathrm{d}\boldsymbol{M}_{N_n^\Lambda(\tau)}+Y_{N_n^\mathrm{M}(\tau-)}^{(n)}\mathrm{d}\boldsymbol{M}_{N_n^\mathrm{M}(\tau)}].
\end{aligned}
\end{equation}
Now we rewrite the system of equations \eqref{36} in the integral form by taking the
expectation and averaging. We have:
\begin{equation}\label{38}
\begin{aligned}
0=&\lim_{T\to\infty}\frac{1}{T}\mathsf{E}\int_{0}^{T}\sum_{j=n+1}^{\infty}\mathsf{I}\{Z(\tau-)=j\}\mathsf{I}\{Y_{N_j^\mathrm{M}(\tau-)+1}^{(j)}=(j-n)\}\mathrm{d}\tau\\
&+\lim_{T\to\infty}\frac{1}{T}\mathsf{E}\int_{0}^{T}\sum_{j=0}^{n-1}\mathsf{I}\{Z(\tau-)=j\}\mathsf{I}\{X_{N_j^\Lambda(\tau-)+1}^{(j)}=(n-j)\}\mathrm{d}\tau\\
&-\lim_{T\to\infty}\frac{1}{T}\mathsf{E}\int_{0}^{T}\mathsf{I}\{Z(\tau-)=n\}\left(X_{N_n^\Lambda(\tau-)}^{(n)}+Y_{N_n^\mathrm{M}(\tau-)}^{(n)}\right)\mathrm{d}\tau\\
&+\lim_{T\to\infty}\frac{1}{T}\mathsf{E}\int_{0}^{T}\sum_{j=n+1}^{\infty}\mathsf{I}\{Z(\tau-)=j\}\mathsf{I}\{Y_{N_j^\mathrm{M}(\tau-)+1}^{(j)}=(j-n)\}\mathrm{d}\boldsymbol{M}_{N_j^M(\tau)}\\
&+\lim_{T\to\infty}\frac{1}{T}\mathsf{E}\int_{0}^{T}\sum_{j=0}^{n-1}\mathsf{I}\{Z(\tau-)=j\}\mathsf{I}\{X_{N_j^\Lambda(\tau-)+1}^{(j)}=(n-j)\}\mathrm{d}\boldsymbol{M}_{N_j^\Lambda(\tau)}\\
&-\lim_{T\to\infty}\frac{1}{T}\mathsf{E}\int_{0}^{T}\mathsf{I}\{Z(\tau-)=n\}[X_{N_n^\Lambda(\tau-)}^{(n)}\mathrm{d}\boldsymbol{M}_{N_n^\Lambda(\tau)}+Y_{N_n^\mathrm{M}(\tau-)}^{(n)}\mathrm{d}\boldsymbol{M}_{N_n^\mathrm{M}(\tau)}]\\
=&J_1+J_2+J_3+J_4+J_5+J_6.
\end{aligned}
\end{equation}
The martingales that appear in the terms $J_4$, $J_5$ and $J_6$ are $\mathscr{F}_\tau$-predictable. Hence the terms $J_4$, $J_5$ and $J_6$ themselves are predictable with respect to $\mathscr{F}=\cup_{\tau>0}\mathscr{F}_{\tau-}$, and
 are equal to zero. For the terms $J_1$ we have
\[
J_1=\lim_{T\to\infty}\frac{1}{T}\int_{0}^{T}\sum_{j=n+1}^{\infty}\mathsf{P}\{Z(\tau-)=j\}\mathsf{P}\{Y_{N_j^\mathrm{M}(\tau-)+1}^{(j)}=j-n\}\mathrm{d}\tau.
\]
Going back to the notation of Example \ref{ex3}, in the limit for $J_1$ we obtain
\[
J_1=\sum_{i=n+1}^{\infty}(i-n)q_{i,n+1}Q_{i}.
\]
Similarly, we obtain
\[
J_2=\sum_{i=0}^{n-1}(n-i)p_{i,n-i}Q_{i},
\]
\[
J_3=\left(\sum_{i=n}^{\infty}(i-n+1)q_{i,n}+\sum_{i=1}^{n}(n-i+1)p_{i,n-i+1}\right)Q_n.
\]
Relations \eqref{e13} follows. The proof of boundary relation \eqref{e13.1} is similar.

From \eqref{e12}--\eqref{e12.1} and \eqref{e13}--\eqref{e13.1} we respectively obtain:
\begin{equation}\label{42}
\begin{aligned}
&P_0+\sum_{n=1}^{\infty}P_n\left(\sum_{i=n}^{\infty}(i-n+1)q_{i,n}+\sum_{i=1}^{n}(n-i+1)p_{i,n-i+1}\right)\\
&=P_1\sum_{i=1}^{\infty}iq_{i,1}+\sum_{n=1}^{\infty}\left(P_{n-1}\sum_{i=0}^{n-1}(n-i)p_{i,n-i}+P_{n+1}\sum_{i=n+1}^{\infty}(i-n)q_{i,n+1}\right),
\end{aligned}
\end{equation}
and
\begin{equation}\label{44}
\begin{aligned}
&Q_0+\sum_{n=1}^{\infty}Q_n\left(\sum_{i=n}^{\infty}(i-n+1)q_{i,n}+\sum_{i=1}^{n}(n-i+1)p_{i,n-i+1}\right)\\
&=Q_1\sum_{i=1}^{\infty}iq_{i,1}+\sum_{n=1}^{\infty}\left(\sum_{i=0}^{n-1}(n-i)p_{i,n-i}Q_{i}+\sum_{i=n+1}^{\infty}(i-n)q_{i,n+1}Q_{i}\right).
\end{aligned}
\end{equation}
Now, consider the expression $\sum_{n=1}^{\infty}\sum_{i=0}^{n-1}(n-i)p_{i,n-i}Q_{i}$. It is not difficult to see that
\begin{equation}\label{46}
\sum_{n=1}^{\infty}\sum_{i=0}^{n-1}(n-i)p_{i,n-i}Q_{i}=\sum_{n=1}^{\infty}Q_{n-1}\sum_{i=0}^{n-1}(n-i)p_{i,n-i}.
\end{equation}
The required result follows by changing the order of the terms in the infinite series, that is legal to do since all the terms of the series are positive.

In addition, for the expression $\sum_{n=1}^{\infty}\sum_{i=n+1}^{\infty}(i-n)q_{i,n+1}Q_{i}$ we have
\begin{equation}\label{48}
\sum_{n=1}^{\infty}\sum_{i=n+1}^{\infty}(i-n)q_{i,n+1}Q_{i}=\sum_{n=1}^{\infty}Q_{n+1}\sum_{i=n+1}^{\infty}(i-n)q_{i,n+1}
\end{equation}
by the same arguments as above.

Relations \eqref{46} and \eqref{48} enable us to represent \eqref{44} as
\begin{equation}\label{50}
\begin{aligned}
&Q_0+\sum_{n=1}^{\infty}Q_n\left(\sum_{i=n}^{\infty}(i-n+1)q_{i,n}+\sum_{i=1}^{n}(n-i+1)p_{i,n-i+1}\right)\\
&=Q_1\sum_{i=1}^{\infty}iq_{i,1}+\sum_{n=1}^{\infty}\left(Q_{n-1}\sum_{i=0}^{n-1}(n-i)p_{i,n-i}+Q_{n+1}\sum_{i=n+1}^{\infty}(i-n)q_{i,n+1}\right).
\end{aligned}
\end{equation}

Hence as before we arrive at the conclusion that $\sum_{n=0}^{\infty}Q_n=1$ if and only if $\sum_{n=0}^{\infty}P_n=1$.
\end{example}

\textit{Derivation of \eqref{e2}--\eqref{e2.2}}.
For the birth-and-death process $\mathscr{B}$, we use the same notation $Z(t)$ for the number of individuals in the population in time $t$, $Z(0)=0$.

For $n\geq2$, the analogue of equation \eqref{12} is
\begin{equation}\label{e6}
\begin{aligned}
&\mathsf{I}\{Z(t)=n\}=\mathsf{I}\{Z(t-)=n+1\}\mathsf{I}\{\mathrm{M}_{n+1}(t)-\mathrm{M}_{n+1}(t-)=1\}\\
&+\mathsf{I}\{Z(t-)=n-2\}\mathsf{I}\{\tilde{\Lambda}_{n-2}(t)-\tilde{\Lambda}_{n-2}(t-)=2\}\\
&+\mathsf{I}\{Z(t-)=n\}\mathsf{I}\{\mathrm{M}_{n}(t)-\mathrm{M}_{n}(t-)=0\}\mathsf{I}\{\tilde{\Lambda}_{n}(t)-\tilde{\Lambda}_{n}(t-)=0\}
\end{aligned}
\end{equation}
that is valid with probability 1,
where $\mathrm{M}_n(t)$ denotes the Poisson process with rate $\mu_n$, $\tilde{\Lambda}_n(t)=2\Lambda_n(t)$, $\Lambda_n(t)$ denotes the Poisson process with rate $\lambda_n/2$.

For $n=0,1$, with probability 1 we have
\begin{equation}\label{e7}
\begin{aligned}
&\mathsf{I}\{Z(t)=1\}=\mathsf{I}\{Z(t-)=2\}\mathsf{I}\{\mathrm{M}_{2}(t)-\mathrm{M}_{2}(t-)=1\}\\
&+\mathsf{I}\{Z(t-)=1\}\mathsf{I}\{\mathrm{M}_{1}(t)-\mathrm{M}_{1}(t-)=0\}\mathsf{I}\{\tilde{\Lambda}_{1}(t)-\tilde{\Lambda}_{1}(t-)=0\},
\end{aligned}
\end{equation}
\begin{equation}\label{e8}
\begin{aligned}
&\mathsf{I}\{Z(t)=0\}=\mathsf{I}\{Z(t-)=1\}\mathsf{I}\{\mathrm{M}_{1}(t)-\mathrm{M}_{1}(t-)=1\}\\
&+\mathsf{I}\{Z(t-)=0\}\mathsf{I}\{\tilde{\Lambda}_{0}(t)-\tilde{\Lambda}_{0}(t-)=0\}.
\end{aligned}
\end{equation}

From \eqref{e6}, \eqref{e7} and \eqref{e8} we arrive at the system of stochastic differential equations, which is an analogue of \eqref{34}:
\begin{equation*}\label{e9}
\begin{aligned}
\mathrm{d}\mathsf{I}\{Z(t)=n\}=&\mathsf{I}\{Z(t-)=n+1\}\mathrm{d}\mathrm{M}_{n+1}(t)\\
&+2\mathsf{I}\{Z(t-)=n-2\}\mathrm{d}\Lambda_{n-2}(t)\\
&-\mathsf{I}\{Z(t-)=n\}(\mathrm{d}\mathrm{M}_{n}(t)+2\mathrm{d}\Lambda_{n}(t)), \ n\geq2,
\end{aligned}
\end{equation*}
\begin{equation*}\label{e10}
\begin{aligned}
\mathrm{d}\mathsf{I}\{Z(t)=1\}=&\mathsf{I}\{Z(t-)=2\}\mathrm{d}\mathrm{M}_{2}(t)\\
&-\mathsf{I}\{Z(t-)=1\}(\mathrm{d}\mathrm{M}_{1}(t)+2\mathrm{d}\Lambda_{1}(t)),
\end{aligned}
\end{equation*}
\begin{equation*}\label{e11}
\mathrm{d}\mathsf{I}\{Z(t)=0\}=\mathsf{I}\{Z(t-)=1\}\mathrm{d}\mathrm{M}_{1}(t)-\mathsf{I}\{Z(t-)=0\}\mathrm{d}\Lambda_{0}(t).
\end{equation*}
The following arguments for the derivation of \eqref{e2}--\eqref{e2.2} are similar to that provided earlier for \eqref{e13}, \eqref{e13.1}.

\smallskip
\textit{Derivation of \eqref{e3}--\eqref{e3.3}} is technically similar to the derivations of the systems of equations provided above.

The three examples considered above show that instead of the original Markov chain model one can study the associated birth-and-death process, if it is only required to find the condition for recurrence or transience rather than exact solution.

\subsection{The final part of the proof}
Let us return to our proof. The proof of the theorem does not distinguish from the proof that has already been provided for the model considered in Example \ref{ex3}. Despite the model considered in Example \ref{ex3} is particular, the derivation of the system of the equations remains the same in general case. Hence, the only thing that is required to do is to consider the associated birth-and-death process and derive the required system of equations.

Denoting the birth-and-death process by $Z(t)$, we have the following system for the state probabilities at time $t+\triangle t$, where $\triangle t$ is a small increment:
\begin{equation}\label{6}
\begin{aligned}
\mathsf{P}\{Z(t+\triangle t)=i\}&=\mathsf{P}\{Z(t+\triangle t)=i~|~A_i^-(\triangle t)\}e_{i+1}^-\triangle t\\
+&\mathsf{P}\{Z(t+\triangle t)=i~|~A_i^+(\triangle t)\}e_{i-1}^+\triangle t\\
+&\mathsf{P}\{Z(t)=i\}\left[1-\triangle t\left(e_i^-+e_i^+\right)\right]+o(\triangle t),
\end{aligned}
\end{equation}
\begin{equation}\label{8}
\begin{aligned}
\mathsf{P}\{Z(t+\triangle t)=0\}&=\mathsf{P}\{Z(t+\triangle t)=0~|~A_1^-(\triangle t)\}e_1^-\triangle t\\
+&\mathsf{P}\{Z(t)=0\}(1-e_0^+\triangle t)+o(\triangle t).
\end{aligned}
\end{equation}

Here in \eqref{6} and \eqref{8}, the meanings of the events $A_i^-(\triangle t)$ and $A_i^+(\triangle t)$  are that during time $(t, t+\triangle t]$ a down-crossing to state $i$ occurs and, respectively, an up-crossing to state $i$ occurs. The meanings of $e_{i}^-$ and $e_{i}^+$, $i\geq1$, are the birth and, respectively, death rates in state $i$; $e_0^+= \sum_{j=1}^\infty jp_{0,j}$. Then the system of equations for the final probabilities of the birth-and-death process follows by the standard way. We have:

\begin{eqnarray*}
P_{n}(e_n^-+e_n^+)&=&P_{n+1}e_{n+1}^-+P_{n-1}e_{n-1}^+, \quad n\geq1,\\
P_0e_0^+&=&P_1{e_1^-}.
\end{eqnarray*}

Thus, the necessary and sufficient condition for recurrence of the Markov chain is
$
\sum_{n=1}^\infty\prod_{i=1}^{n}e_i^-/e_i^+=\infty.
$
The theorem is proved.

\begin{remark}
The reduction to the system of equations that describes birth-and-death processes enables us to adapt the other known results on recurrence and transience given in \cite{A1, A2} as well.
\end{remark}

\section{Examples}\label{S3}
In this section, we provide elementary but non-trivial examples of recurrent and transient Markov chains with an infinite set of states.

1. Consider the following matrix of transition probabilities
\begin{equation*}
\left(\begin{matrix}0 &1 &0 &0 &0 &0 &0 &\cdots\\
\smallskip
\frac{7}{12} &0 &\frac{3}{12} &\frac{2}{12} &0 &0 &0 &\cdots\\
\smallskip
0 &\frac{7}{12} &0 &\frac{3}{12} &\frac{2}{12} &0 &0 &\cdots\\
\smallskip
0 &0 &\frac{7}{12} &0 &\frac{3}{12} &\frac{2}{12} &0 &\cdots\\
\smallskip
\vdots &\vdots &\vdots &\ddots &\ddots &\ddots &\ddots &\vdots
\end{matrix}\right)
\end{equation*}

The associated transition probability matrix of the birth-and-death type is
\begin{equation*}
\left(\begin{matrix}0 &1 &0 &0 &0 &0  &\cdots\\
\smallskip
\frac{7}{10} &0 &\frac{3}{10} &0 &0  &0 &\cdots\\
\smallskip
0 &\frac{1}{2} &0 &\frac{1}{2} &0  &0 &\cdots\\
\smallskip
0 &0 &\frac{1}{2} &0 &\frac{1}{2}  &0 &\cdots\\
\smallskip
\vdots &\vdots &\vdots &\ddots &\ddots &\ddots &\vdots
\end{matrix}\right)=
\left(\begin{matrix}0 &1 &0 &0 &0 &0  &\cdots\\
\smallskip
q_1 &0 &p_1 &0 &0 &0  &\cdots\\
\smallskip
0 &q_2 &0 &p_2 &0 &0  &\cdots\\
\smallskip
0 &0 &q_3 &0 &p_3 &0  &\cdots\\
\smallskip
\vdots &\vdots &\vdots &\ddots &\ddots &\ddots &\vdots
\end{matrix}\right).
\end{equation*}
With $q_n=p_n=1/2$ for $n\geq2$, the birth-and-death process is recurrent. Hence, the above Markov chain is recurrent.

2. In the following example that is similar to the example above, the orders of fractions $3/12$ and $2/12$ in the matrix are changed:
\begin{equation*}
\left(\begin{matrix}0 &1 &0 &0 &0 &0 &0 &\cdots\\
\smallskip
\frac{7}{12} &0 &\frac{2}{12} &\frac{3}{12} &0 &0 &0 &\cdots\\
\smallskip
0 &\frac{7}{12} &0 &\frac{2}{12} &\frac{3}{12} &0 &0 &\cdots\\
\smallskip
0 &0 &\frac{7}{12} &0 &\frac{2}{12} &\frac{3}{12} &0 &\cdots\\
\smallskip
\vdots &\vdots &\vdots &\ddots &\ddots &\ddots &\ddots &\vdots
\end{matrix}\right)
\end{equation*}
In this case, the Markov chain is transient.

Indeed, the associated transition probability matrix for the birth-and-death process now is

\begin{equation*}
\left(\begin{matrix}0 &1 &0 &0 &0 &0  &\cdots\\
\smallskip
\frac{7}{9} &0 &\frac{2}{9} &0 &0  &0 &\cdots\\
\smallskip
0 &\frac{7}{15} &0 &\frac{8}{15} &0  &0 &\cdots\\
\smallskip
0 &0 &\frac{7}{15} &0 &\frac{8}{15}  &0 &\cdots\\
\smallskip
\vdots &\vdots &\vdots &\ddots &\ddots &\ddots &\vdots
\end{matrix}\right)=
\left(\begin{matrix}0 &1 &0 &0 &0 &0  &\cdots\\
\smallskip
q_1 &0 &p_1 &0 &0 &0  &\cdots\\
\smallskip
0 &q_2 &0 &p_2 &0 &0  &\cdots\\
\smallskip
0 &0 &q_3 &0 &p_3 &0  &\cdots\\
\smallskip
\vdots &\vdots &\vdots &\ddots &\ddots &\ddots &\vdots
\end{matrix}\right).
\end{equation*}
Now we have $p_n>q_n$ for all $n\geq2$, and hence the transience of the birth-and-death process and the related Markov chain follows.

\section{Counterexamples}\label{S4}
The aim of this section is to show that the basic assumptions of the theorem that the states of the Markov chain form connected domain is important, and if this assumption is not satisfied the statement of the theorem fails.

Indeed, let us consider the following matrix of transition probabilities:
\begin{equation*}\label{18}
\left(\begin{matrix}0 &1 &0 &0 &0 &0 &0 &0 &0 &\ldots\\
\epsilon &0 &1-\epsilon &0 &0 &0 &0 &0 &0 &\ldots\\
\frac{1}{3}-\epsilon &\epsilon &0 &\epsilon &\frac{2}{3}-\epsilon &0 &0 &0 &0 &\ldots\\
0 &\frac{1}{2}+\epsilon &0 &0 &0  &\frac{1}{2}-\epsilon &0 &0 &0 &\ldots\\
0 &0 &\frac{1}{3}-\epsilon &\epsilon &0 &\epsilon &\frac{2}{3}-\epsilon &0 &0 &\ldots\\
0 &0 &0 &\frac{1}{2}+\epsilon &0 &0 &0  &\frac{1}{2}-\epsilon &0 &\ldots\\
\vdots &\vdots &\vdots &\vdots &\ddots &\ddots &\ddots &\ddots &\ddots &\vdots\end{matrix}\right),
\end{equation*}
where $\epsilon>0$ is a small fixed parameter. Note that the considered matrix is irreducible. However the Markov chain does not form connected domain. According to this matrix of transition probabilities, the typical excursion of the Markov chain is as follows. From state $0$ it visits state $1$. Then, it visits state $2$ or otherwise goes back to the initial state. Appearing once in an even state (in our case it is state $2$), the further excursion through even states is as follows:
from any state with even order number $2n$, the chain visits the state $2n+2$ with probability that is slightly less than $2/3$, or the state $2n-2$ with probability that is slightly less than $1/3$, and there is small probability to visit a state with an odd order number. An excursion over even states can be long, but once an odd state will be visited. Once visiting an odd state (distinct from $1$), the chain finally will return to state $1$, since the probability of down-crossing (which is with step 2) is higher than the such one of up-crossing. The excursion that starts from state $1$ can be repeated many times, but once the chain will return to the initial state. Our conclusion is that the Markov chain is recurrent. On the other hand, it is readily seen that in this case $\sum_{n=1}^\infty\prod_{i=1}^{n}{e_i^-}/{e_i^+}<\infty$. So, the statement of the theorem cannot be applied and fails.

The example above demonstrated the case of recurrent Markov chain, for which the condition of the theorem is not satisfied. It is not difficult to demonstrate a similar example for a transient Markov chain for which the condition of the theorem is not satisfied either. Here is such an example
\begin{equation*}\label{20}
\left(\begin{matrix}0 &1 &0 &0 &0 &0 &0 &0 &0 &\ldots\\
1-\epsilon &0 &\epsilon &0 &0 &0 &0 &0 &0 &\ldots\\
\frac{2}{3}-\epsilon &\epsilon &0 &\epsilon &\frac{1}{3}-\epsilon &0 &0 &0 &0 &\ldots\\
0 &\frac{1}{2}-\epsilon &0 &0 &0  &\frac{1}{2}+\epsilon &0 &0 &0 &\ldots\\
0 &0 &\frac{2}{3}-\epsilon &\epsilon &0 &\epsilon &\frac{1}{3}-\epsilon &0 &0 &\ldots\\
0 &0 &0 &\frac{1}{2}-\epsilon &0 &0 &0  &\frac{1}{2}+\epsilon &0 &\ldots\\
\vdots &\vdots &\vdots &\vdots &\ddots &\ddots &\ddots &\ddots &\ddots &\vdots\end{matrix}\right),
\end{equation*}
which we leave without comments because of similarity.

\section{Discussion}\label{S5} In this paper, we established a new result for recurrence and transience for relatively wide class of Markov chains. We showed, that if the basic condition of the theorem is not satisfied, then its statement fails. The counterexamples that built in Section \ref{S3} are important. They shows that even if a Markov chain is irreducible, that is space states forms a single communicating class, there may be separate subclasses that qualitatively affect on the behavior of the Markov chain. In our case, one subclass includes all even states, another one all odd states. In these subclasses the behavior of the Markov chain demonstrated in the first example  is as follows. In the even states, the state order number change rapidly up, whereas in the odd states the state order number changes slowly down. Because of specified connection between the two subclasses of states, it turns out that state $1$ is accessed with probability $1$, from which the Markov chain eventually returns to the initial state with probability $1$. In the second example, the behavior is opposite. In even states the state numbers change rapidly down, whereas in odd states the state numbers change slowly up increasing to infinity.

In our opinion, the following future study might be interesting.

1. The main result of this paper assumes that the matrix of transition probabilities form connected domain. However, for many transition probability matrices that do not form connected domain the results remain true, and the proof provided in this paper does not use this class of matrices explicitly. \textit{Can the class of transition probability matrices be extended?}

2. One of the conditions of Theorem \ref{Th1} requires that the diagonal elements of the matrix satisfy the condition: \textit{$1-p_{i,i}>\epsilon$ for all $i\geq0$, where $\epsilon$ is some positive value.} Can this condition be relaxed in some way? Say, one can assume that $p_{i,i}$ tends to $1$ as $i\to\infty$, but at some slow rate that will enable us to claim that the main result still holds true.

\subsection*{Data availability statement}
Data sharing not applicable to this article as no datasets were generated or analysed during the current study.

\subsection*{Disclosures and declarations}
No conflict of interest was reported by the author.

\subsection*{Authorship clarified}
The author, who conducted this work, has no relation to any institution.

\bibliographystyle{amsplain}

\end{document}